\documentclass[12pt]{article}

\usepackage{epsfig, graphicx, psfrag, subfigure}

\usepackage{amsfonts}
\usepackage{amsmath}
\usepackage{amssymb}
\usepackage{latexsym}
\usepackage{amscd, amsthm}
\usepackage{bbm}
\usepackage{url}

\newcommand{\R}{\mathbb{R}}
\newcommand{\E}{\mathbb{E}}

\newcommand{\Prob}{\mathbb{P}}

\newcommand{\Cov}{\mathrm{Cov}}
\newcommand{\Vol}{\mathrm{Vol \,}}
\newcommand{\Id}{\mathrm{Id}}
\newcommand{\Var}{\mathrm{Var}}
\newcommand{\conv}{\mathrm{Conv \,}}

\newcommand{\eps}{\varepsilon}
\newcommand{\te}{\theta}
\newcommand{\al}{\alpha}

\newtheorem{theorem}{Theorem}
\newtheorem{corollary}{Corollary}
\newtheorem{Conjecture}{Conjecture}

\begin{document}
\title{Concentration phenomena in high dimensional geometry.}
\author{Olivier Gu\'edon}
\newcommand\address{\noindent\leavevmode\noindent
Olivier Gu\'{e}don, \\
Universit\'{e} Paris-Est \\
Laboratoire d'Analyse et Math\'{e}matiques Appliqu\'ees (UMR 8050). \\
UPEMLV, F-77454, Marne-la-Vall\'ee, France  \\
\texttt{\small e-mail: olivier.guedon@univ-mlv.fr}
}

\date{}
\maketitle

\begin{abstract}The purpose of this note is to present several aspects of concentration phenomena in high dimensional geometry. At the heart of the study  is a geometric analysis point of view coming from the theory of high dimensional convex bodies.  The topic has a broad audience going from algorithmic convex geometry to random matrices. We have tried to emphasize different problems relating these areas of research.  Another connected area is the study of probability in Banach spaces where some concentration phenomena are related with good comparisons between the weak and the strong moments of a random vector.  
\end{abstract}

\section{Convex geometry and log-concave measures}
A function
$f : \R^n \to \R_+$ is said to be log-concave if $\forall x, y \in \R^n, \forall \theta \in [0,1],$
$$
f((1-\theta) x + \theta y) \ge f(x)^{1- \theta} f(y)^{\theta}
$$
Define a measure $\mu$ with a log-concave density $f \in L_1^{\mathrm{loc}}$, the Pr\'ekopa-Leindler inequality 
 \cite{MR0404557, MR2199372} implies that it satisfies: for every compact sets $A , B \subset \R^n$, for every $\theta \in [0,1]$,
\begin{equation}
\label{mu:logconcave}
\mu((1-\theta) A + \theta B) \ge \mu(A)^{1- \theta} \mu(B)^{\theta}.
\end{equation}
A measure satisfying $(\ref{mu:logconcave})$ is said to be log-concave. A complete characterization of log-concave measures is well known. It has been done during the sixties and seventies and it is related to the work of \cite{MR0056669, MR0315079, MR0404557, MR2199372,  MR0404559}, see also the surveys \cite{MR588074, MR2167203}. In \cite{MR0404559} it is proved that a measure is log-concave if and only if it is absolutely continuous with respect to the Lebesgue measure on the affine subspace generated by its convex support, with log-concave locally integrable density. Classical examples are the case of product of exponential distributions, $f(x) = \frac{1}{2^n} \exp(-|x|_1)$, the Gaussian measure, $f(x) = \frac{1}{(2\pi)^{n/2}} \exp(-|x|_2^2/2)$, the uniform measure on a convex body, $f(x) = 1_K(x)$.  Moreover, it is well known \cite{MR0241584, MR0404557, MR2199372, MR0388475} that the class of log-concave measures is stable under convolution and linear transformations. It is important to notice that the class of uniform distribution on a convex body is stable under linear transformation but not under convolution. This is one among several reasons why it is preferable to work with log-concave measures.
\subsection{The hyperplane conjecture or the slicing problem.}
It is one of the famous conjecture in high dimensional convex geometry. 
\begin{Conjecture} {\sc (The hyperplane conjecture)}
\label{slicing}
There exist a constant $C>0$ such that
for every $n$ and every convex body $K \subset \R^n$ of volume 1 and barycenter at the origin, there is a direction $\theta$ such that ${{\rm Vol} \, (K \cap \theta^{\perp}) \ge C}$.
\end{Conjecture}
Several other formulations are known. For example, it is equivalent to ask if for 
every convex bodies  $K_1$ and $K_2$  with barycenter at the origin such that for every 
$\theta \in S^{n-1}$
${\rm Vol} \, (K_1 \cap \theta^{\perp}) \le {\rm Vol} \, (K_2 \cap \theta^{\perp})$ 
then  ${\rm Vol} \, (K_1) \le \, C \,  {\rm Vol} \, (K_2)$ ? It is worth noticed that the constant $C$ can not be 1 in high dimension. Indeed, replacing $C$ by 1 in the conclusion lead to the Busemann-Petty problem which is known to be true in dimension $n \le 4$ but false  in dimension $n \ge 5$ \cite{MR1689343}. 
We refer to \cite{MR1008717} for a more detailed presentation of the slicing problem. 
For a convex body $K$, define $L_K$ by 
\begin{equation}
\label{def:LK}
n \, L_K^2 = \min_{T \in SL_n(\R)} \frac{1}{({\rm Vol} \, K)^{1+\frac{2}{n}}} 
\int_K |Tx|_2^2 \, dx
\end{equation}
where the minimum is taken over all affine transformations preserving the volume.
Observe that $L_K$ is invariant under affine transformation. From a result of Hensley \cite{MR572315},  the hyperplane conjecture is equivalent to the following question. Does there exist a constant $C$ such that for every dimension $n$ and convex body $K \subset \R^n$, $L_K \le C \ ?$ The number $L_K$ is called the isotropic  constant of the convex body $K$. 
The minimum in $(\ref{def:LK})$ is attained when the ellipsoid $\cal E$ is the inertia matrix associated to $K$, centered at the barycenter of $K$. Equivalently, it is attained for $B_2^n$ when $K$ is in isotropic position, that is : 
$$
\frac{1}{{\rm Vol} \, K}  \int_K x dx = 0, \quad  \hbox{ and } \quad \frac{1}{{\rm Vol} \, K} \int_K x_i x_j \,dx = \delta_{i,j}. 
$$ 
In isotropic position, we have $L_K = \frac{1}{({\rm Vol} \, K)^{\frac{1}{n}}}$. 
It is also possible to define the isotropic constant of a log-concave function $f : \R^n \to \R_+$. Let $X$ be the random vector in $\R^n$ with density of probability $f$, define 
$$
L_f = L_X = \left( \det(\Cov X) \right)^{2/n} \left(f(\E X) \right)^{1/n}
$$
It is not difficult to check that it is invariant under linear transformation, $L_X = L_{T(X)}$ for any $T \in GL_n(\R)$ so that we can assume that $X$ is isotropic. This means that 
$$
\int f(x)dx = 1, \quad \int x f(x) = 0 \quad  \hbox{ and } \quad  \int x_i x_j \, f(x) dx = \delta_{i,j}
$$
and that $L_X = f(0)^{1/n}$.
A  theorem of Ball \cite{MR932007} asserts that the hyperplane conjecture is equivalent to the uniform boundedness of the isotropic constant of log-concave random vectors.  
\begin{theorem}[Ball \cite{MR932007}]
\label{th:Ball}
There exists a constant $C$ such that 
\[
\sup_{n, K} L_K \le \sup_{n,f} L_f \le C \sup_{n,K} L_K
\]
where the suprema are taken with respect to every dimension $n$, every (isotropic) convex bodies and every (isotropic) log-concave functions. 
\end{theorem}
The left hand side of the inequality is obvious. There are more important ingredients to prove the right hand side inequality.  In particular, it requires to define a convex body from a log-concave measure and to keep a good control of the isotropic constant. Let $f : \R^n \to \R_+$ be a log-concave function and define a family of set $\left(K_p(f)\right)_{p>0}$ by
\begin{equation}
\label{eq:Kp}
K_p(f) = \left\{
x \in \R^n, \int_0^{+\infty} t^{p-1} f(tx) dt \ge \frac{f(0)}{p}
\right\}
\end{equation}
A main step in the proof of Ball's theorem above is to prove that $K_p(f)$ is a convex set for every $p>0$. Moreover, there are some good relations between $L_f$ and the isotropic constant of $K_{n+1}(f)$.  The best known bound today is due to Klartag \cite{MR2276540} who proved that for any log-concave function $f : \R^n \to \R_+$, $L_f \le c n^{1/4}$.
\subsection{Computing the volume of a convex body}
In algorithmic convex geometry \cite{MR1261419}, an important question is to produce algorithm which may compute the volume of a convex body. It is another  illustration  of a high dimensional phenomena in convex geometry. To start to understand the question, we need to be more precise. We assume that the convex body $K$ is given by a separation oracle. A separation oracle has the property that if you give a point $x \in \R^n$ then it answers either that $x \in K$ or that $x \notin K$ and it gives a separating hyperplane between $x$ and $K$. There are powerful negative results in this direction, see \cite{MR866364, MR911186}. These results show that the volume of a convex body can not be even approximated by any  deterministic polynomial algorithm. Be aware that  in fixed dimension $n$, these problems can easily be solved in polynomial time. However, the degree of those polynomials estimating the running time increases fast with $n$. What I know about this negative result is that it is based  on an important theorem of high dimensional convex geometry due independently to Carl-Pajor \cite{MR969241}, Gluskin \cite{MR945901} and Barany-F\"uredi \cite{MR911186}. It evaluates the volume of the absolute convex hull of $N$ points $u_1, \ldots, u_N$ from the unit sphere in $\R^n$,  
\[
\left(
\frac{\Vol  \conv \{\pm u_1, \ldots, \pm u_N  \}}{\Vol B_2^n} 
\right)^{1/n} \leq C \ \sqrt {\frac{\log\left(1+ \frac{N}{n}\right)}{n}} .
\]
In the topic of algorithmic convex geometry, a breakthrough has been done by Dyer-Frieze-Kannan \cite{MR1095916}. They proved that the situation changed drastically if we allow randomization and they gave the first polynomial randomized algorithm to calculate the volume of a convex body. It exists today a vast literature on the subject and we refer to the paper of Kannan, Lov\'asz and Simonovits \cite{MR1608200} and to \cite{MR1261419, MR2178341, MR2853824} for a better presentation of this topic.
In the randomized situation,  we are given two numbers $\eps$ and $\eta$ in $(0,1)$, a convex body $K$ and a separation oracle. Then there exists a randomized algorithm returning a non negative number $\zeta$ such that 
\[
(1- \eps) \zeta < \Vol K < (1+ \eps) \zeta.
\]
with probability at least $1-\eta$. What is random in this situation is the number of oracle calls  and they get (for example) a bound on its expected value. In \cite{MR1608200}, the algorithm uses 
\[
O\left( \frac{n^5}{\eps^2} \left( \log \frac{1}{\eps} \right)^{3} \left( \log \frac{1}{\eta} \right) \log^5 n \right)
\]
random call to the separation oracle 
which was a significant improvement of the previous known results. The strategy from \cite{MR1608200} may be described in two steps and each of them led to  interesting problems in high dimensional convex geometry. 

The first step is a rounding procedure in order to insure that the convex body and the Euclidean structure are correctly related. 
It consists to finding  a position such that $B_2^n \subset K \subset d \, B_2^n$ 
where $d$ depends polynomially in the dimension $n$. A theorem of John \cite{MR0030135} asserts that there exists an affine transformation $T$ such that 
$B_2^n$ is the ellipsoid of maximal volume contained in $T(K)$ and in which case, 
$B_2^n \subset T(K)  \subset n B_2^n$ (the distance is reduced to $\sqrt n$ in the case of symmetric convex bodies). However, there is no 
known
randomized polynomial algorithm that can return the affine transformation $T$ such that $B_2^n$ is close to the ellipsoid of maximal volume contained in $K$. A classical procedure was to use the ellipsoid algorithm (see \cite{MR1261419}) which achieves $d = O(n^{3/2})$. 
The main idea is to consider the inertia ellipsoid $\cal E$ associated to $K$.  It is possible to approximate it by a randomized polynomial time algorithm. Moreover, the geometric distance between $K$ and $\cal E$ is known to be of the order of $n$ and it is not difficult to prove that a big part of the volume of $K$, $(1-\eps) \Vol K$, is contained in $\left(\log \frac{1}{\eps} \right) \sqrt n \, {\cal E}$. The question that attracted a lot of attention in high dimensional convex geometry was to find a good random way to approximate the inertia ellipsoid. This has been solved recently by Adamczak, Litvak, Pajor, Tomczak-Jaegermann in \cite{MR2601042} and we will give a more detailed description of the problem in the paragraph \ref{sub:inertia}.

The second step consists of computing the volume of a convex body which is in a nearly isotropic position. They apply  a multiphase Monte-Carlo algorithm.  
They need good bounds on the mixing time of the random walk and this is based on good estimate about an isoperimetric inequality. We will not describe the notions of local conductance here and refer to \cite{MR1318794, MR1608200}. We will focus on an isoperimetric problem for convex bodies in isotropic position. The question of describing the "almost" extremal sets in the isoperimetric inequality is still an open problem, known today as the KLS conjecture \cite{MR1318794}. We will give a more detailed description of the problem in the paragraph \ref{sub:isoperimetry}.

\subsection{Approximation of the inertia matrix}
\label{sub:inertia}
Let $X$ be a random vector uniformly distributed on a convex body in $\R^n$. The inertia matrix is given by $\E X \otimes X$. The simplest procedure to approximate the inertia matrix is to understand how many sample is needed to approximate it. Given $\eps \in (0,1)$, the question is to give an estimate of the smallest number $N$ such that 
\[
\left\| \frac{1}{N} \sum_{i=1}^N X_i \otimes X_i - \E X \otimes X \right\| \le \eps \| \E X \otimes X \|
\]
where $\| \cdot \|$ denotes the operator norm from $\ell_2^n$ to $\ell_2^n$. Since the procedure is random, we would like to have such result with some fixed positive probability, $1-\eta$.  It is clear that without loss of generality, we can assume that the random vector is isotropic, that is 
$\E X = 0$ and $\E X \otimes X = \Id$. In terms of random processes, the question becomes: evaluate $N$ such that with the highest possible probability, 
\begin{equation}
\label{eq:approximationidentity}
\sup_{y \in S^{n-1}} \left|  \frac{1}{N} \sum_{i=1}^N \langle X_j, y \rangle^2 - 1 \right|
\le
\eps .
\end{equation}
Using the language of random matrices, this is nothing else that evaluating $N$ such that all the singular values of the random matrix $A$, with rows $X_1 / \sqrt N, \ldots, X_N / \sqrt N$, belong to the interval $[1-\eps, 1+ \eps]$. In \cite{MR1608200} Kannan, Lov\'asz and Simonovits  proved that if $N$ is of the order of 
$\frac{n^2}{\eps^2  \eta^2}$ then $(\ref{eq:approximationidentity})$ holds true with probability larger than $1-\eta$.  Shortly after, Bourgain \cite{MR1665576} improved this estimate to $N \approx \frac{n \log^3 n}{\eps^2 \eta^2} $. During more than fifteen years, several groups of people \cite{MR1694526, MR1800254, MR2190337, MR2276533, MR2304336} proposed different strategies to improve this result. A  breakthrough has been done by Adamczak, Litvak, Pajor and Tomczak-Jaegermann \cite{MR2601042, MR2769907}  who proved that if 
$N \approx \frac{n}{\eps^2}$ then $(\ref{eq:approximationidentity})$ holds true with probability at least $1 - e^{-c \sqrt n}$. The main achievement in the result is that $N$ is taken of the order of the ambient dimension $n$ and that the probability of the event is not only large but increases extremely fast regarding to the dimension $n$. It is now of interest to understand what other random vector than log-concave probability distribution satisfy such type of result. It has been recently investigated by Srivastava and Vershynin \cite{SV} and generalized to the case of random matrices in \cite{Youssef}. 

\subsection{Almost extremal sets in the isoperimetric inequality}
\label{sub:isoperimetry}
The isoperimetric problem for convex bodies is the following. Let $K$ be a convex body in $\R^n$ and $\mu$ be the uniform measure on $K$. Let $S$ be a subset of $K$ and define the boundary measure of $S$ as 

\vskip0.2cm
$
\includegraphics[height=1in,width=1.9in]{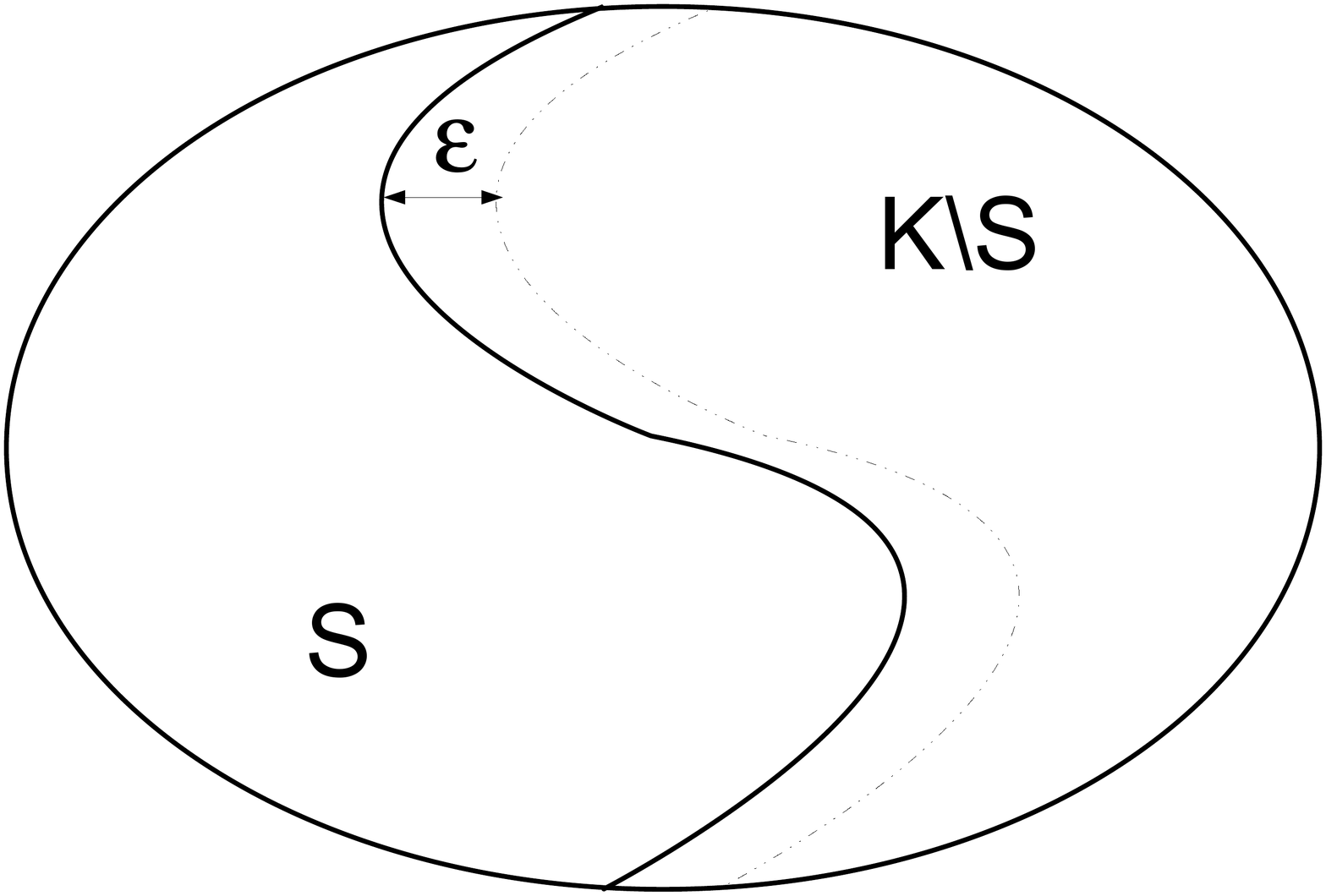}
$

\vskip-2cm
\hskip 5.95cm 
$\displaystyle{
\mu^+(S) = \liminf_{\eps \to 0} \frac{\mu(S + \eps B_2^n) - \mu(S)}{\eps}.
}$

\vskip1cm
\noindent
This definition is also valid for any measure with log-concave density on $\R^n$. The question is to evaluate the largest possible $h$ such that 
\begin{equation}
\label{eq:isoperimetrie}
 \forall \ S \subset K,  \ \mu^+(S) \ge \ h \ \mu(S) (1-\mu(S))
\end{equation}
Without any assumption on the measure, you can easily imagine a situation where $h$ may be as close to $0$ as you wish. In our situation, we made  the assumption that the measure is isotropic and log-concave. This avoid a lot of non regular situation. Kannan, Lov\'asz and Simonovits (the same group of persons in an earlier paper !) \cite{MR1318794} considered this problem and conjectured that up to a universal constant in the inequality $(\ref{eq:isoperimetrie})$, the worth set $S$ should be a half space of the same measure than $S$. 
\begin{Conjecture} {\sc (The KLS conjecture)}
\label{KLS}
There exists $c>0$ such that for any dimension $n$ and any isotropic log-concave probability on $\R^n$,
\[
 \forall \ S \subset \R^n,  \ \mu^+(S) \ge \ c \ \mu(S) (1-\mu(S))
\]
\end{Conjecture}
This is supported by the fact that in the Gaussian setting, it is known since the work of Sudakov, Tsirelson \cite{MR0365680} and independently Borell \cite{MR0399402}  that the half spaces are the exact solutions  of the isoperimetric problem in the Gauss space: 
\[
\min_{\gamma_n(S) = \alpha}  \gamma_n(S + \eps B_2^n) = \gamma_n(H + \eps B_2^n) = \int_{-\infty}^{a + \eps} e^{-t^2/2}  \frac{dt}{\sqrt{2\pi}}
\]
where $d\gamma_n(x)= e^{-|x|_2^2/2} dx / (2\pi)^{n/2}$ and $H$ is a half space of Gaussian measure 
\[
\alpha = \int_{-\infty}^{a} e^{-t^2/2}  \frac{dt}{\sqrt{2\pi}}.
\] 
The inequality $(\ref{eq:isoperimetrie})$ is called a Cheeger type inequality and $h$ is usually referred as the Cheeger's constant of the measure $\mu$. Kannan, Lov\'asz and Simonovits proved that if $\mu$ is a log-concave isotropic probability measure on $\R^n$, then 
\[
h \ge \frac{c}{\E |X|_2} \approx \frac{1}{\sqrt n}
\]
where $X$ is a random vector distributed according to $\mu$.
The argument is based on a localization technique introduced in \cite{MR901393, MR1238906}. The localization method was further developed in \cite{MR2060645, MR2249622}. And Bobkov \cite{MR2347041} improved this result to 
\begin{equation}
\label{eq:Bobkov}
h \ge \frac{c}{\left( \Var |X|^2\right)^{1/4}}
\end{equation}
During the last decade, this question has been much investigated. However, it remains an open question. Very few positive results are known. It has been proved only for some classes of convex bodies like the unit balls of $\ell_p^n$ \cite{MR2446328, MR2449135} and a weaker form is proved for random Gaussian polytopes in \cite{MR2875755}. It is also known from the work of Buser \cite{MR683635} and Ledoux \cite{MR1186991} that in the case of log-concave probability, the Cheeger constant  is related to the best constant in the Poincar\'e inequality. Let $X$ be the random vector distributed according to $\mu$, let $D_2$ be the largest constant  such that for every regular function $F : \R^n \to \R$,
\begin{equation}
\label{eq:PoincareL2}
D_2 \ \Var F(X) \le \E |\nabla F(X)|_2^2
\end{equation}
then $h^2 \approx D_2$. More surprinsigly, Milman \cite{MR2507637} proved that $h^2 \approx D_{\infty}$ where $D_{\infty}$ is the largest constant  such that for every 1-Lipschitz function $F : \R^n \to \R$,
\[
\Var F(X) \le \frac{1}{D_\infty}.
\]
Both inequalities are easy consequences of Conjecture \ref{KLS}. The difficult part of the proof of these results concern the reverse statement. A recent paper of Gozlan, Roberto and Samson \cite{GRS} completes  also the picture of the different equivalent formulations of the question. 

We refer to Chapter 3 in \cite{MR1849347} and to \cite{MR901393} for a more detailed description of the links between the Poincar\'e inequality and concentration of measure and we just emphasize on the fact that Conjecture \ref{KLS} implies a very strong  concentration inequality of the Euclidean norm.
\begin{Conjecture} {\sc (The thin shell conjecture)}
\label{Concentration}
There exists $c>0$ such that 
for any log-concave isotropic probability on $\R^n$, for any $t>0$,
\[
\Prob \left( \bigg| |X|_2 - \sqrt n \bigg| \ge t \sqrt n \right) \le 2 e^{- c\, t \sqrt n}
\]
\end{Conjecture}
It is the purpose of the next section to describe our knowledge about this question and the different aspects of this problem in high dimensional convex geometry.
\section{The thin shell concentration and a Berry-Esseen type theorem for convex bodies}
The classical Berry-Esseen bounds of the central limit theorem asserts that for every independent random variables $x_1, \ldots, x_n$ such that $\E x_i = 0$, $\E x_i^2 = 1$ and $\E x_i^3 = \tau$, we have 
\[
\forall \te \in S^{n-1}, \quad
\sup_{t \in \R} \left| \Prob\left(\sum_{i=1}^n \te_i x_i \le t \right)
-
\int_{-\infty}^t e^{-u^2/2} \frac{du}{\sqrt{2 \pi}}\right|
\le 
\tau |\theta|_4^2 .
\]
Observe that if $\te = (1/\sqrt n, \ldots, 1/\sqrt n)$, we have $ |\theta|_4^2 = 1/ \sqrt n$ which gives a very good rate of convergence. At the end of the nineties, Ball (in several talks and in \cite{MR1997580}) and independently Brehm and Voigt in \cite{MR1801435} asked about a possible generalization of the Berry-Esseen Theorem for convex bodies. The question can be stated as follows:  for every isotropic convex body $K \subset \R^n$, does there exist a direction $\te \in S^{n-1}$ such that
\begin{equation}
\label{eq:CLT}
\sup_{t \in \R} \left| \Prob\left(\sum_{i=1}^n \te_i x_i \le t \right)
-
\int_{-\infty}^t e^{-u^2/2} \frac{du}{\sqrt{2 \pi}}\right|
\le 
\al_n
\end{equation}
with $\lim_{+\infty} \al_n = 0$, where $\Prob$ is the uniform distribution on $K$ ? The same question can be asked for every  isotropic probability $\Prob$ with log-concave density on $\R^n$.  

A very good presentation of the problem is done in \cite{MR1997580}. Indeed the authors propose a satisfactory way to solve the problem and this generated an intense activity during the last decade. Anttila, Ball and Perssinaki \cite{MR1997580} proved that if there exists a sequence $(\eps_n)_{n \ge 1}$ such that $\lim_{+ \infty} \eps_n = 0$ and 
\begin{equation}
\label{eq:ABP}
\Prob\left( \left| \frac{|X|_2}{\sqrt n} - 1 \right| \ge \eps_n \right) \le \eps_n
\end{equation}
for every isotropic probability with log-concave density on $\R^n$, then the Berry-Esseen Theorem for convex bodies holds true. More precisely, they proved that if $(\ref{eq:ABP})$ is true for an isotropic probability $\Prob$ uniformly distributed on a convex body (or with log-concave density) then there exists a set $A$ of directions $\theta \in S^{n-1}$ of extremely large probability such that $(\ref{eq:CLT})$ holds true.  They also give an affirmative answer to $(\ref{eq:ABP})$ for some classes of uniform measures like the uniform measures on  the unit ball of $\ell_p^n$ for all $p\ge 1$. From a more probabilistic point of view, the principle of studying a thin shell estimate to get a Central Limit Theorem appeared already in \cite{MR1440135, MR751274}. 
More precise formulation of Berry-Esseen type theorem for isotropic random vector satisfying $(\ref{eq:ABP})$ have been established by Bobkov \cite{MR1959791}. 
The fact that under the isotropicity condition, for most of the directions $\theta \in S^{n-1}$, the random variables $\langle X, \theta \rangle$ have a common behavior was already observed in \cite{MR517198}.

A complete solution for the Central Limit Problem for convex bodies is given by Klartag in \cite{MR2285748}, while shortly after, an independent proof of $(\ref{eq:ABP})$ was given in \cite{MR2349721}. In these papers, the rates for $\al_n$ and $\eps_n$ are weak, only of the order of some negative power of $\log n$. In fact, the inequality $(\ref{eq:ABP})$ is  equivalent to  
\begin{equation}
\label{eq:Var}
\lim_{n \to \infty} \frac{\Var |X|_2}{n} = 0
\end{equation}
where $X$ is distributed according to the log-concave measure $\Prob$ and this is known to be true in full generality \cite{MR2285748, MR2349721}.
As it was emphasized by Bobkov and Koldobsky \cite{MR2083387},  if we look at the particular case of $F$ being $|\cdot |_2^2$ in the conjecture of 
Kannan, Lov\'asz and Simonovits, we deduce from inequality 
$(\ref{eq:PoincareL2})$ the following conjecture.
\begin{Conjecture} {\sc (The Variance conjecture)}
\label{ThinShell}
There exists a constant $C$ such that for every isotropic log-concave random vector $X$, we have
$\Var |X|_2^2 \le C \E |X|_2^2$ or equivalently 
\[
\left(\E |X|_2^4\right)^{1/4} \le \left( 1 + \frac{C}{n}\right)\left(\E |X|_2^2 \right)^{1/2}. 
\]
\end{Conjecture}
It seems to be a very natural conjecture. Denoting $X = (x_1, \ldots, x_n)$, it concerns the behavior of the random variable $x_1^2 + \ldots + x_n^2$, where here, the classical hypotheses on the entries of $X$, independence and variance 1,  are replaced by the isotropy and log-concavity of the vector. 
Taking $F$ as the Euclidean norm in $(\ref{eq:PoincareL2})$, this gives another simple conjecture
\begin{Conjecture} {\sc (The Variance conjecture - bis)}
\label{ThinShellbis}
There exists a constant $C$ such that for every isotropic log-concave random vector $X$, 
\[
\Var |X|_2 \le C
\]
\end{Conjecture}
It is clear that these conjectures are much stronger than inequality $(\ref{eq:Var})$. 

Let us come back to Conjecture \ref{Concentration}. The stated concentration inequality contains informations for all $t > 0$ and the discussion is different for $t > 10$, for $t$ close to 1 and for $t < 1/10$. Let us start to discuss about the large deviation, $t > 10$. From Borell's lemma \cite{MR0388475} (see also \cite{MR856576} Appendix 3), it is easy to deduce that for $t > 10$, $\Prob ( |X|_2 \ge t \sqrt n ) \le e^{-ct}$. But it was proved in \cite{MR2083388} that in the  case of uniform measure on an unconditional isotropic convex body, this can be improved to $e^{-ct \sqrt n}$. In this topic, a breakthrough has been done by Paouris \cite{MR2276533}. Shortly after, he proved in \cite{MR2833584} a small ball inequality which corresponds to the case of $t$ being close to $1$. 
\begin{theorem}[Paouris \cite{MR2276533, MR2833584}]
\label{th:Paouris}
Let $\Prob$ be an isotropic probability on $\R^n$ with log-concave density. 
\\
Then for every $t \ge 10$, 
\[
\Prob( |X|_2 \ge t \sqrt n ) \le e^{-ct \sqrt n}.
\]
And for every $\eps \in (0,1/10)$,
\[
\Prob( |X|_2 \le \eps \sqrt n) \le (c \eps)^{\sqrt n}
\]
\end{theorem}
An important ingredient in his proof is the study of the volume of the sets $K_p(f)$ defined in $(\ref{eq:Kp})$, the Keith Ball's bodies introduced in \cite{MR932007} to prove Theorem \ref{th:Ball}.
\\
Once the thin shell estimate $(\ref{eq:ABP})$ has been proved by Klartag \cite{MR2285748}, see also \cite{MR2349721}, it was natural to study the rate of convergence.
The results have been improved to polynomial estimates in the dimension $n$ by Klartag \cite{MR2311626}, by Fleury \cite{MR2652173} and lastly by E. Milman and myself  \cite{MR2846382}. As of today, the 
best known result is the following 
\begin{theorem}[Gu\'edon-Milman \cite{MR2846382}]
\label{th:GM}
\[
\forall t \ge 0, \quad
\Prob \left( \big| |X|_2 -\sqrt{n} \big | \geq t \sqrt{n}  \right) \leq C \exp(-c  \sqrt n \ \min(t^3,t)) 
\]
\end{theorem}
From this result we deduced that 
$\Var |X|_2^2 \le C \, n^{5/3}$ and using $(\ref{eq:Bobkov})$, this gives a general bound on the Cheeger's constant of an isotropic log-concave measure  
$h \ge c \, n^{-5/12}$. A new relation between the thin shell estimate and the Cheeger's constant has been recently developed by Eldan \cite{Eldan2012} where he deduced $h \ge c \, n^{-1/3} (\log n)^{-1/2}$ from Theorem \ref{th:GM}.

It is time to draw a picture of a high dimensional convex body. This picture was popularized by Vitali Milman. Observe in particular that it is important to draw the convex set as a star shape body with a lot of points very far from the origin and lot of points very close to the origin. Think to the example of the cube whose $2^n$ vertices are at Euclidean distance  $\sqrt n$ of the origin and whose $2n$ middle of faces are at Euclidean distance $1$.

\vskip-0.6cm

\includegraphics[height=4in,width=4in]{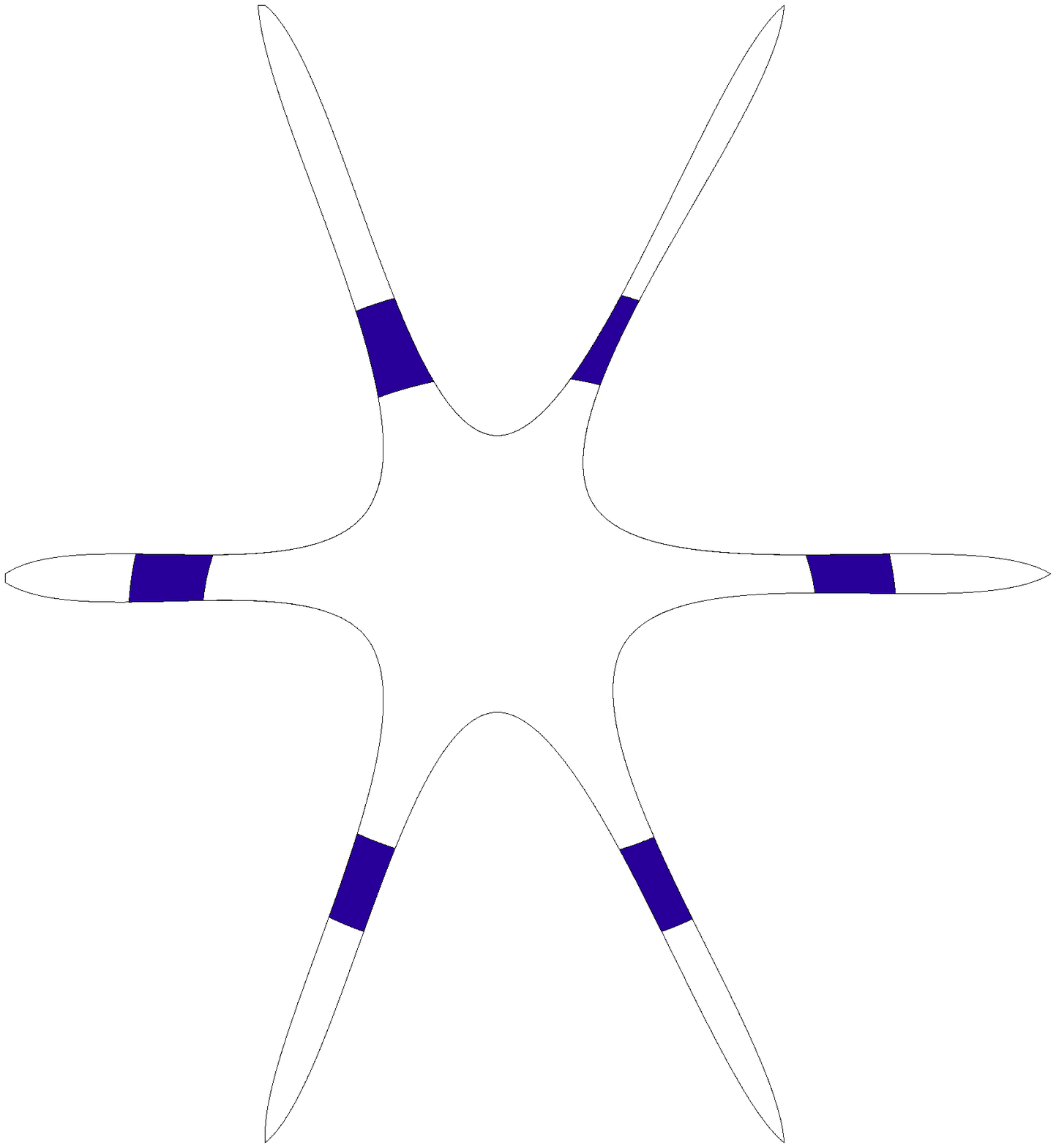}

\vskip-8cm
\hskip 7.5cm 
A high dimensional isotropic convex body. 

\vskip 2cm
\hskip 7cm 
The volume is concentrated in the thin shell  of

\hskip 7cm 
radius $n^{1/2}$ and width $n^{1/2-1/6}$ drawn in blue.

\vskip 2cm
\noindent
Since the volume of the convex body is highly concentrated in the blue part of the picture, it maybe explains the difficulty of computing it with some algorithm because you need to generate a random walk in several disconnected parts.

The concentration of the mass of a log-concave measure in a  Euclidean ball or Euclidean shell is usually  understood via  the study of the $L_p$ norms of the random variable $|X|_2$. A new proof of the first  part of Theorem \ref{th:Paouris} has been recently given in \cite{ALLOPT}. It is worth noticing that the general result can be expressed without the isotropic hypothesis. It gives a deep relation of the strong moments of $|X|_2$ with its weak moments defined by $\sigma_p^p(X) = \sup_{|z|_2 \le 1}  \E |\langle z, X \rangle |^p .$
\begin{theorem}[\cite{MR2276533, ALLOPT}]
\label{th:ALLOPT}
There exist $c, C > 0$ such that for every $\log$-concave random vector $X$, we have  
\begin{equation}
\label{eq:ALLOPT}
 \forall p \ge 1, 
 \quad 
\left( \E |X|_2^p\right)^{1/p}
\le 
C \  \E |X|_2 + c \, \sigma_p(X) 
\end{equation}
where $\sigma_p(X) = \sup_{|z|_2 \le 1} \left( \E | \langle z, X \rangle |^p \right)^{1/p}.$
\end{theorem}
Observe that in isotropic position, $\E|X|_2 \le (\E|X|_2^2)^{1/2} = \sqrt n $. 
By Borell's lemma \cite{MR0388475} (see also \cite{MR856576} Apppendix 3), we know that 
\[
\forall p \ge 1, \quad \left( \E |\langle z, X \rangle |^p \right)^{1/p}  \le C \, p \,
\left( \E \langle z, X \rangle^2 \right)^{1/2} = C \, p \, |z|_2
\]
Hence by Theorem \ref{th:ALLOPT}, for all $p \ge 1$, 
$\left( \E |X|_2^p\right)^{1/p}
\le 
C \sqrt n  + c  p$.
Take $p = t \sqrt n$ and  Markov inequality gives
\[
\forall t \ge 1, \quad
\Prob\left( |X|_2 \ge t \, \sqrt n\right) \le e^{-c \, t \, \sqrt n}.
\]
This is exactly the first part of Theorem \ref{th:Paouris}. In \cite{MR2846382} we also investigated the $L_p$ norms of $|X|_2$. Looking at Conjecture \ref{ThinShell} and at the statement of the previous Theorem,   we should be able to prove the following weaker conjecture.
\begin{Conjecture} {\sc (The weak thin shell conjecture)}
\label{WeakMoment}
There exists $c > 0$ such that for every $\log$-concave random vector $X$, we have 
\[
 \forall p \ge 1,
 \quad 
\left( \E |X|_2^p\right)^{1/p}
\le 
 \E |X|_2 + c \, \sigma_p(X) 
\]
\end{Conjecture}
This would prove that in isotropic position 
\[
 \forall p \ge 1,
 \quad 
\left( \E |X|_2^p\right)^{1/p}
\le 
 \E |X|_2 \left(1  + c  \, \frac{p}{\sqrt n} \right).
\]
Recently, Eldan and Klartag \cite{MR2858465} established a surprising connection between thin shell estimates and the slicing problem. They proved that if 
Conjecture \ref{ThinShell} holds true then Conjecture \ref{slicing} will be also true. It is also known \cite{BallNguyen} that for an individual 
log-concave distribution $\mu$, the isotropic constant of $\mu$ is bounded by a function of the Cheeger's constant of $\mu$. 
Finally, Eldan \cite{Eldan2012} proved that Conjecture \ref{ThinShell} 
implies Conjecture \ref{KLS} up to polylogarithmic term. 

To conclude this paragraph, I would like to advertise that Conjecture \ref{Concentration} is known to be true when the measure is uniformly distributed on the unit ball of a generalized Orlicz space \cite{MR2667700} and that Conjecture \ref{ThinShell} is valid when the measure is uniformly distributed on an unconditional convex body \cite{MR2520120}.

\section{Weak and strong moments of a random vector}
In this section, we will go away from the framework defined by log-concave distributions.
We have seen that the thin shell conjectures for log-concave measures are closely related to the comparison of the strong moments of the Euclidean norm with its weak moments. There are more probabilistic questions  in this direction. It is natural to consider for which family of random vectors $X$, we have  for any norm
$\| \cdot \|$
\[
(\E\|X\|^p)^{1/p} \le C \, \E\|X\| + c \, \sup_{\|z\|_{\star} \le 1} (\E\langle z , X\rangle^p)^{1/p}
\]
where $\| \cdot \|_{\star}$ is the dual norm of $\| \cdot \|$ and $C, c >0$ are numerical constants. It is known to be true with $C =1$ for Gaussian \cite{MR0458556} or Rademacher 
(see \cite{MR2814399} Theorem 4.7) random vectors series.  We refer to \cite{MR2449135} and \cite{MR2918093} where such questions are discussed. In the area of log-concave measures, Lata\l a \cite{MR2918093} asked the following.
\begin{Conjecture} {\sc (Weak and strong moments conjecture)}
\label{Latala}
There exist $C, c > 0$ such that for any $\log$-concave random vector $X$ and any norm $\| \cdot \|$ we have  
\[
(\E\|X\|^p)^{1/p} \le C \, \E\|X\| + c \, \sup_{\|z\|_{\star} \le 1} (\E\langle z , X\rangle^p)^{1/p}
\]
where $\| \cdot \|_{\star}$ is the dual norm of $\| \cdot \|$.
\end{Conjecture}
Observe that Theorem \ref{th:ALLOPT} tells that this is true for the Euclidean norm. From now on, we will present a few of the results from a joint work with Adamczak, Lata\l a, Litvak, Pajor and Tomczak-Jaegermann \cite{AGLLOPT}. We introduce a new class of random vectors. 
Let $p>0$, $m=\lceil p\rceil$, and $\lambda\geq 1$. We say that a
random vector $X$ in a Banach space $E$  satisfies the assumption $H(p,\lambda)$ 
if for every linear mapping $A:E\to\R^m$ such that  $Y=AX$ is 
non-degenerate
there exists a gauge $\|\cdot\|$ on $\R^m$ such that
$\E\|Y\|< \infty$ and 
\[
(\E\|Y\|^p)^{1/p}\leq \lambda\,\E \|Y\|.
\]
Be aware that in the definition, we chose $m=\lceil p\rceil$. And observe that any $m$-dimensional norm may be approximated by a set $\{\varphi_i\}$ of $3^m$ linear forms which means that for any norm in $\R^m$, we have 
\[
(\E\|Y\|^p)^{1/p}\ \approx \left( \E  \sup_{i = 1, \ldots, 3^m} | \varphi_i(Y)|^p \right)^{1/p} \approx \sup_{\|\varphi\|_{\star} \le 1 } \left( \E   | \varphi(Y)|^p \right)^{1/p}.
\]
We say that a random vector $X \in E$ is $\psi_2$ with constant $\psi$ if it satisfies
\[
\forall p \ge 2, \forall \phi  \in E^*, \quad  \left( \E   | \varphi(Y)|^p \right)^{1/p} \le \psi \sqrt p \left( \E   | \varphi(Y)|^2 \right)^{1/2}.
\]
From the previous remark, we can easily deduce that $\psi_2$ random vectors satisfy the hypothesis $H(p, C \psi^2)$ for every $p$ with a universal constant $C$.  This implies that any Gaussian or Rademacher random series of vectors satisfy the hypothesis $H(p, C)$ for every $p$. 
The first main result from \cite{AGLLOPT} is
\begin{theorem}
\label{th:main1}
Let $p > 0$ and $\lambda \ge 1$. If a random vector $X$ satisfies $H(p, \lambda)$ then 
$$
(\E |X|_2^p )^{1/p} \le c \, ( \lambda \E |X|_2 + \sigma_p(X) )
$$
where $c$ is a universal constant.
\end{theorem}

\noindent {\bf Proof.}
A very rough idea of the proof is the following. Let 
$X$ be the random vector in $E$ satisfying $H(p, \lambda)$, $m = \lceil p \rceil$, $\lambda \ge 1$. There are three important steps in the proof.
\\
By Gaussian Concentration \cite{MR0458556}, we know that for  a standard Gaussian vector $G$
\begin{equation}
\label{eq:concentration}
(\E_G \E_X | \langle G, X \rangle |^p)^{1/p}
\le 
\E_G (\E_X | \langle G, X \rangle |^p)^{1/p} + c \, \sqrt p \ \sigma_p(X)
\end{equation}
We define a new norm on $E$ by 
\[
\| z \| = ( \E_X | \langle z, X \rangle |^p)^{1/p}.
\]
It is nothing else than the dual norm of the classical $Z_p$ bodies \cite{MR2276533}. 
By Gordon min-max theorem \cite{MR800188}, we know  that if $A : (E, \|\cdot\|) \to \R^m \sim {\cal N}(0, \Id)$ is a standard Gaussian matrix, we have
\begin{equation}
\label{eq:Gordon}
\E_G  ( \E_X | \langle G, X \rangle |^p)^{1/p}
\le \E_A \min_{|z|_2 = 1}  ( \E_X | \langle z, AX \rangle |^p)^{1/p} + c \, \sqrt p \ \sigma_p(X).
\end{equation}
The key property $H(p, \lambda)$ allows to prove the following nice geometric Lemma
\begin{equation}
\label{eq:geometric}
\min_{|z|_2 = 1}  ( \E_X | \langle z, AX \rangle |^p)^{1/p}
\le \lambda \ \E_X |AX|_2.
\end{equation}
To conclude, it remains to glue the argument. Observe that 
\[
 (\E|X|_2^p)^{1/p} \approx \frac{1}{\sqrt p} (\E_G \E_X | \langle G, X \rangle |^p)^{1/p}.
\]
By $(\ref{eq:concentration})$ and $(\ref{eq:Gordon})$, we get
\[
(\E|X|_2^p)^{1/p} \approx \frac{1}{\sqrt p} (\E_G \E_X | \langle G, X \rangle |^p)^{1/p}
\lesssim
\frac{1}{\sqrt p}  \E_A \min_{|z|_2 = 1}  ( \E_X | \langle z, AX \rangle |^p)^{1/p}  +  \ \sigma_p(X).
\]
From $(\ref{eq:geometric})$, we conclude that
\begin{align*}
(\E|X|_2^p)^{1/p} & \approx \frac{1}{\sqrt p} (\E_G \E_X | \langle G, X \rangle |^p)^{1/p}
\leq
\frac{1}{\sqrt p}  \E_A \min_{|z|_2 = 1}  ( \E_X | \langle z, AX \rangle |^p)^{1/p}  +  \ \sigma_p(X)
\\
& \lesssim 
\frac{1}{\sqrt p} \E_A \, \lambda \ \E_X |AX|_2 +  \ \sigma_p(X)
\lesssim
\lambda \ \E |X|_2 + \sigma_p(X).
\end{align*}
This is the announced result.
\qed

Of course, we need to understand which random vectors satisfy the hypothesis $H(p, \lambda)$ with $p$ as large as possible and $\lambda$ being a constant. It is not only the case for Gaussian, Rademacher random series of vectors, or $\psi_2$ random vectors. It is satisfied by log-concave random vectors (as expected) and more generally by $s$-concave random vectors for negative $s$. The class of $s$-concave random vectors is very general. For example, the Cauchy distributions belong to it. It has been widely studied in the seventies and there is a serie of papers of Borell \cite{MR0404559, MR0388475} who characterized it. 
\\
Let $s < 1/n$. A probability Borel measure $\mu$ on $\R^n$ is called $s$-concave if
for every compact sets $A , B \subset \R^n$, for every $\theta \in [0,1]$,
\[
\mu((1-\theta) A + \theta B) \ge ((1- \theta)\mu(A)^{s} + \theta \mu(B)^{s})^{1/s}
\]  
whenever $\mu(A) \mu(B) > 0$.
\\
For $s=0$, this corresponds to log-concave measures.
\\
Borell \cite{MR0404559, MR0388475} characterized the class of $s$-concave measures and proved that any $s$-concave probability is supported on some convex subset of an affine subspace where it has a density. Assuming the measure has full dimensional support, the density is a $\gamma$-concave function
with $\frac{1}{\gamma} = \frac{1}{s} - n$.
We will  concentrate on the case of $s$-concave random vectors with $s<0$. 
Denote $s = -1/r$. Restating the characterization of Borell, we get that
when the support generates the whole space, the $s$-concave measure has a density $g$ of the form
\[
g = f^{-\beta} \quad \hbox{ with } \quad \beta = n+r
\]
and $f$ is a positive convex function on $\R^n$. 
A classical example is to define $g$ from a norm $\| \cdot \|$ on $\R^n$ :
\[
g(x) = c (1 + \|x\|)^{-n-r}, r > 0.
\]
The class of $s$-concave measures satisfy several important properties. It is decreasing in $s$ which gives for example that 
any log-concave probability measure is $(-1/r)$-concave for any $r>0$. 
The linear image of a $(-1/r)$-concave vector is also $(-1/r)$-concave.
And the Euclidean norm of a $(-1/r)$-concave random vector has moments of order $0<p<r$.
Our main second result \cite{AGLLOPT} is 
\begin{theorem}
Let $r \ge 2$ and $X$ be a $(-1/r)$-concave random vector. Then for every 
$0<p<r/2$, $X$ satisfies the assumption $H(p,C)$, $C$ being a universal constant.
\end{theorem}
We refer to \cite{AGLLOPT} for its proof. With Theorem \ref{th:main1}, we get
\begin{theorem}
Let $r \ge 2$ and $X$ be a $(-1/r)$-concave random vector. Then for every 
$0<p<r/2$, 
$$
(\E|X|_2^p)^{1/p} \le C (\E |X|_2 + \sigma_p(X)).
$$
\end{theorem}
It is not difficult to deduce nice concentration properties of the random variable $|X|_2$ when $X$ is an isotropic $(-1/r)$-concave random vector in $\R^n$. 
\begin{corollary}
Let $r \ge 2$ and $X$ be an isotropic $(-1/r)$-concave random vector in $\R^n$. Then for every $t > 0$,
\[
\Prob\big( |X|_2 > t \sqrt n \big)
\le 
\left(
\frac{c \max(1, r / \sqrt n)}{t}
\right)^{r/2}
\]
where $c$ is a universal constant.
\end{corollary}
This concentration property is central in the work of 
Srivastava and Vershynin \cite{SV} about approximation of the covariance matrix. This is why we can go back to the problem studied in section \ref{sub:inertia}
about the approximation of the inertia matrix. 
\begin{corollary}
Let $r \ge \log n$ and $X$ be a $(-1/r)$-concave isotropic random vector. Let $X_1, \ldots, X_N$ be independent copies of $X$. Then for every $\eps \in (0,1)$ and every $N \ge C(\eps) n$, one has
$$
\E \left \| \frac{1}{N} \sum_{i=1}^N X_i \otimes X_i^{} - I \right \| \le \eps.
$$
\end{corollary}
Of course, we stated only a few of the results from \cite{AGLLOPT} and refer to it for more general and precise formulations. In particular, we also got some interesting small ball properties in the spirit of the results of Paouris \cite{MR2833584}.  We also asked if some thin shell concentration could be proved for $-1/r$-concave random vector and this was done very recently in \cite{FGP}.

{\footnotesize
\bibliographystyle{acm}
\bibliography{GuedonBibTex}

\def\cprime{$'$} \def\cprime{$'$} \def\cprime{$'$}
\begin{thebibliography}{10}

\bibitem{AGLLOPT}
{\sc Adamczak, R., Gu\'edon, O., Lata\l{a}, R., Litvak, A.~E., Oleszkiewicz,
  K., Pajor, A., and Tomczak-Jaegermann, N.}
\newblock Moment estimates for convex measures.
\newblock {\em Electron. J. Probab. 17\/} (2012), 1--19.

\bibitem{ALLOPT}
{\sc Adamczak, R., Lata\l{a}, R., Litvak, A.~E., Oleszkiewicz, K., Pajor, A.,
  and Tomczak-Jaegermann, N.}
\newblock A short proof of paouris' inequality.
\newblock {\em Can. Math. Bul.\/} (to appear).

\bibitem{MR2601042}
{\sc Adamczak, R., Litvak, A.~E., Pajor, A., and Tomczak-Jaegermann, N.}
\newblock Quantitative estimates of the convergence of the empirical covariance
  matrix in log-concave ensembles.
\newblock {\em J. Amer. Math. Soc. 23}, 2 (2010), 535--561.

\bibitem{MR2769907}
{\sc Adamczak, R., Litvak, A.~E., Pajor, A., and Tomczak-Jaegermann, N.}
\newblock Sharp bounds on the rate of convergence of the empirical covariance
  matrix.
\newblock {\em C. R. Math. Acad. Sci. Paris 349}, 3-4 (2011), 195--200.

\bibitem{MR1997580}
{\sc Anttila, M., Ball, K., and Perissinaki, I.}
\newblock The central limit problem for convex bodies.
\newblock {\em Trans. Amer. Math. Soc. 355}, 12 (2003), 4723--4735
  (electronic).

\bibitem{MR932007}
{\sc Ball, K.}
\newblock Logarithmically concave functions and sections of convex sets in
  {${\bf R}^n$}.
\newblock {\em Studia Math. 88}, 1 (1988), 69--84.

\bibitem{BallNguyen}
{\sc Ball, K., and Nguyen, V.~H.}
\newblock Entropy jumps for random vectors with log-concave density and
  spectral gap.
\newblock {\em http://arxiv.org/abs/1206.5098\/} (preprint).

\bibitem{MR911186}
{\sc B{\'a}r{\'a}ny, I., and F{\"u}redi, Z.}
\newblock Computing the volume is difficult.
\newblock {\em Discrete Comput. Geom. 2}, 4 (1987), 319--326.

\bibitem{MR1959791}
{\sc Bobkov, S.~G.}
\newblock On concentration of distributions of random weighted sums.
\newblock {\em Ann. Probab. 31}, 1 (2003), 195--215.

\bibitem{MR2347041}
{\sc Bobkov, S.~G.}
\newblock On isoperimetric constants for log-concave probability distributions.
\newblock In {\em Geometric aspects of functional analysis}, vol.~1910 of {\em
  Lecture Notes in Math.} Springer, Berlin, 2007, pp.~81--88.

\bibitem{MR2083387}
{\sc Bobkov, S.~G., and Koldobsky, A.}
\newblock On the central limit property of convex bodies.
\newblock In {\em Geometric aspects of functional analysis}, vol.~1807 of {\em
  Lecture Notes in Math.} Springer, Berlin, 2003, pp.~44--52.

\bibitem{MR2083388}
{\sc Bobkov, S.~G., and Nazarov, F.~L.}
\newblock On convex bodies and log-concave probability measures with
  unconditional basis.
\newblock In {\em Geometric aspects of functional analysis}, vol.~1807 of {\em
  Lecture Notes in Math.} Springer, Berlin, 2003, pp.~53--69.

\bibitem{MR0388475}
{\sc Borell, C.}
\newblock Convex measures on locally convex spaces.
\newblock {\em Ark. Mat. 12\/} (1974), 239--252.

\bibitem{MR0399402}
{\sc Borell, C.}
\newblock The {B}runn-{M}inkowski inequality in {G}auss space.
\newblock {\em Invent. Math. 30}, 2 (1975), 207--216.

\bibitem{MR0404559}
{\sc Borell, C.}
\newblock Convex set functions in {$d$}-space.
\newblock {\em Period. Math. Hungar. 6}, 2 (1975), 111--136.

\bibitem{MR1665576}
{\sc Bourgain, J.}
\newblock Random points in isotropic convex sets.
\newblock In {\em Convex geometric analysis ({B}erkeley, {CA}, 1996)}, vol.~34
  of {\em Math. Sci. Res. Inst. Publ.} Cambridge Univ. Press, Cambridge, 1999,
  pp.~53--58.

\bibitem{MR1801435}
{\sc Brehm, U., and Voigt, J.}
\newblock Asymptotics of cross sections for convex bodies.
\newblock {\em Beitr\"age Algebra Geom. 41}, 2 (2000), 437--454.

\bibitem{MR683635}
{\sc Buser, P.}
\newblock A note on the isoperimetric constant.
\newblock {\em Ann. Sci. \'Ecole Norm. Sup. (4) 15}, 2 (1982), 213--230.

\bibitem{MR969241}
{\sc Carl, B., and Pajor, A.}
\newblock Gel\cprime fand numbers of operators with values in a {H}ilbert
  space.
\newblock {\em Invent. Math. 94}, 3 (1988), 479--504.

\bibitem{MR0458556}
{\sc Cirel{\cprime}son, B.~S., Ibragimov, I.~A., and Sudakov, V.~N.}
\newblock Norms of {G}aussian sample functions.
\newblock In {\em Proceedings of the {T}hird {J}apan-{USSR} {S}ymposium on
  {P}robability {T}heory ({T}ashkent, 1975)\/} (Berlin, 1976), Springer,
  pp.~20--41. Lecture Notes in Math., Vol. 550.

\bibitem{MR588074}
{\sc Das~Gupta, S.}
\newblock Brunn-{M}inkowski inequality and its aftermath.
\newblock {\em J. Multivariate Anal. 10}, 3 (1980), 296--318.

\bibitem{MR0241584}
{\sc Davidovi{\v{c}}, J.~S., Korenbljum, B.~I., and Hacet, B.~I.}
\newblock A certain property of logarithmically concave functions.
\newblock {\em Dokl. Akad. Nauk SSSR 185\/} (1969), 1215--1218.

\bibitem{MR751274}
{\sc Diaconis, P., and Freedman, D.}
\newblock Asymptotics of graphical projection pursuit.
\newblock {\em Ann. Statist. 12}, 3 (1984), 793--815.

\bibitem{MR1095916}
{\sc Dyer, M., Frieze, A., and Kannan, R.}
\newblock A random polynomial-time algorithm for approximating the volume of
  convex bodies.
\newblock {\em J. Assoc. Comput. Mach. 38}, 1 (1991), 1--17.

\bibitem{Eldan2012}
{\sc Eldan, R.}
\newblock Thin shell implies spectral gap up to polylog via a stochastic
  localization scheme.
\newblock {\em http://arxiv.org/abs/1203.0893\/} (preprint).

\bibitem{MR2858465}
{\sc Eldan, R., and Klartag, B.}
\newblock Approximately {G}aussian marginals and the hyperplane conjecture.
\newblock In {\em Concentration, functional inequalities and isoperimetry},
  vol.~545 of {\em Contemp. Math.} Amer. Math. Soc., Providence, RI, 2011,
  pp.~55--68.

\bibitem{MR866364}
{\sc Elekes, G.}
\newblock A geometric inequality and the complexity of computing volume.
\newblock {\em Discrete Comput. Geom. 1}, 4 (1986), 289--292.

\bibitem{MR2667700}
{\sc Fleury, B.}
\newblock Between {P}aouris concentration inequality and variance conjecture.
\newblock {\em Ann. Inst. Henri Poincar\'e Probab. Stat. 46}, 2 (2010),
  299--312.

\bibitem{MR2652173}
{\sc Fleury, B.}
\newblock Concentration in a thin {E}uclidean shell for log-concave measures.
\newblock {\em J. Funct. Anal. 259}, 4 (2010), 832--841.

\bibitem{MR2875755}
{\sc Fleury, B.}
\newblock Poincar\'e inequality in mean value for {G}aussian polytopes.
\newblock {\em Probab. Theory Related Fields 152}, 1-2 (2012), 141--178.

\bibitem{MR2349721}
{\sc Fleury, B., Gu{\'e}don, O., and Paouris, G.}
\newblock A stability result for mean width of {$L_p$}-centroid bodies.
\newblock {\em Adv. Math. 214}, 2 (2007), 865--877.

\bibitem{MR2060645}
{\sc Fradelizi, M., and Gu{\'e}don, O.}
\newblock The extreme points of subsets of {$s$}-concave probabilities and a
  geometric localization theorem.
\newblock {\em Discrete Comput. Geom. 31}, 2 (2004), 327--335.

\bibitem{MR2249622}
{\sc Fradelizi, M., and Gu{\'e}don, O.}
\newblock A generalized localization theorem and geometric inequalities for
  convex bodies.
\newblock {\em Adv. Math. 204}, 2 (2006), 509--529.

\bibitem{FGP}
{\sc Fradelizi, M., Gu\'edon, O., and Pajor, A.}
\newblock Spherical thin-shell concentration for convex measures.

\bibitem{MR1689343}
{\sc Gardner, R.~J., Koldobsky, A., and Schlumprecht, T.}
\newblock An analytic solution to the {B}usemann-{P}etty problem on sections of
  convex bodies.
\newblock {\em Ann. of Math. (2) 149}, 2 (1999), 691--703.

\bibitem{MR2190337}
{\sc Giannopoulos, A., Hartzoulaki, M., and Tsolomitis, A.}
\newblock Random points in isotropic unconditional convex bodies.
\newblock {\em J. London Math. Soc. (2) 72}, 3 (2005), 779--798.

\bibitem{MR1800254}
{\sc Giannopoulos, A.~A., and Milman, V.~D.}
\newblock Concentration property on probability spaces.
\newblock {\em Adv. Math. 156}, 1 (2000), 77--106.

\bibitem{MR945901}
{\sc Gluskin, E.~D.}
\newblock Extremal properties of orthogonal parallelepipeds and their
  applications to the geometry of {B}anach spaces.
\newblock {\em Mat. Sb. (N.S.) 136(178)}, 1 (1988), 85--96.

\bibitem{MR800188}
{\sc Gordon, Y.}
\newblock Some inequalities for {G}aussian processes and applications.
\newblock {\em Israel J. Math. 50}, 4 (1985), 265--289.

\bibitem{GRS}
{\sc Gozlan, N., Roberto, C., and Samson, P.-M.}
\newblock From dimension free concentration to poincar\'e inequality.
\newblock {\em http://arxiv.org/abs/1305.4331\/} (preprint).

\bibitem{MR901393}
{\sc Gromov, M., and Milman, V.~D.}
\newblock Generalization of the spherical isoperimetric inequality to uniformly
  convex {B}anach spaces.
\newblock {\em Compositio Math. 62}, 3 (1987), 263--282.

\bibitem{MR1261419}
{\sc Gr{\"o}tschel, M., Lov{\'a}sz, L., and Schrijver, A.}
\newblock {\em Geometric algorithms and combinatorial optimization},
  second~ed., vol.~2 of {\em Algorithms and Combinatorics}.
\newblock Springer-Verlag, Berlin, 1993.

\bibitem{MR2846382}
{\sc Gu{\'e}don, O., and Milman, E.}
\newblock Interpolating thin-shell and sharp large-deviation estimates for
  isotropic log-concave measures.
\newblock {\em Geom. Funct. Anal. 21}, 5 (2011), 1043--1068.

\bibitem{MR2304336}
{\sc Gu{\'e}don, O., and Rudelson, M.}
\newblock {$L_p$}-moments of random vectors via majorizing measures.
\newblock {\em Adv. Math. 208}, 2 (2007), 798--823.

\bibitem{MR572315}
{\sc Hensley, D.}
\newblock Slicing convex bodies---bounds for slice area in terms of the body's
  covariance.
\newblock {\em Proc. Amer. Math. Soc. 79}, 4 (1980), 619--625.

\bibitem{MR0056669}
{\sc Henstock, R., and Macbeath, A.~M.}
\newblock On the measure of sum-sets. {I}. {T}he theorems of {B}runn,
  {M}inkowski, and {L}usternik.
\newblock {\em Proc. London Math. Soc. (3) 3\/} (1953), 182--194.

\bibitem{MR0030135}
{\sc John, F.}
\newblock Extremum problems with inequalities as subsidiary conditions.
\newblock In {\em Studies and {E}ssays {P}resented to {R}. {C}ourant on his
  60th {B}irthday, {J}anuary 8, 1948}. Interscience Publishers, Inc., New York,
  N. Y., 1948, pp.~187--204.

\bibitem{MR1318794}
{\sc Kannan, R., Lov{\'a}sz, L., and Simonovits, M.}
\newblock Isoperimetric problems for convex bodies and a localization lemma.
\newblock {\em Discrete Comput. Geom. 13}, 3-4 (1995), 541--559.

\bibitem{MR1608200}
{\sc Kannan, R., Lov{\'a}sz, L., and Simonovits, M.}
\newblock Random walks and an {$O^*(n^5)$} volume algorithm for convex bodies.
\newblock {\em Random Structures Algorithms 11}, 1 (1997), 1--50.

\bibitem{MR2276540}
{\sc Klartag, B.}
\newblock On convex perturbations with a bounded isotropic constant.
\newblock {\em Geom. Funct. Anal. 16}, 6 (2006), 1274--1290.

\bibitem{MR2285748}
{\sc Klartag, B.}
\newblock A central limit theorem for convex sets.
\newblock {\em Invent. Math. 168}, 1 (2007), 91--131.

\bibitem{MR2311626}
{\sc Klartag, B.}
\newblock Power-law estimates for the central limit theorem for convex sets.
\newblock {\em J. Funct. Anal. 245}, 1 (2007), 284--310.

\bibitem{MR2520120}
{\sc Klartag, B.}
\newblock A {B}erry-{E}sseen type inequality for convex bodies with an
  unconditional basis.
\newblock {\em Probab. Theory Related Fields 145}, 1-2 (2009), 1--33.

\bibitem{MR2918093}
{\sc Lata{\l}a, R.}
\newblock Weak and strong moments of random vectors.
\newblock In {\em Marcinkiewicz centenary volume}, vol.~95 of {\em Banach
  Center Publ.} Polish Acad. Sci. Inst. Math., Warsaw, 2011, pp.~115--121.

\bibitem{MR2449135}
{\sc Lata{\l}a, R., and Wojtaszczyk, J.~O.}
\newblock On the infimum convolution inequality.
\newblock {\em Studia Math. 189}, 2 (2008), 147--187.

\bibitem{MR1186991}
{\sc Ledoux, M.}
\newblock A simple analytic proof of an inequality by {P}. {B}user.
\newblock {\em Proc. Amer. Math. Soc. 121}, 3 (1994), 951--959.

\bibitem{MR1849347}
{\sc Ledoux, M.}
\newblock {\em The concentration of measure phenomenon}, vol.~89 of {\em
  Mathematical Surveys and Monographs}.
\newblock American Mathematical Society, Providence, RI, 2001.

\bibitem{MR2814399}
{\sc Ledoux, M., and Talagrand, M.}
\newblock {\em Probability in {B}anach spaces}.
\newblock Classics in Mathematics. Springer-Verlag, Berlin, 2011.
\newblock Isoperimetry and processes, Reprint of the 1991 edition.

\bibitem{MR2199372}
{\sc Leindler, L.}
\newblock On a certain converse of {H}\"older's inequality. {II}.
\newblock {\em Acta Sci. Math. (Szeged) 33}, 3-4 (1972), 217--223.

\bibitem{MR1238906}
{\sc Lov{\'a}sz, L., and Simonovits, M.}
\newblock Random walks in a convex body and an improved volume algorithm.
\newblock {\em Random Structures Algorithms 4}, 4 (1993), 359--412.

\bibitem{MR2167203}
{\sc Maurey, B.}
\newblock In\'egalit\'e de {B}runn-{M}inkowski-{L}usternik, et autres
  in\'egalit\'es g\'eom\'etriques et fonctionnelles.
\newblock {\em Ast\'erisque}, 299 (2005), Exp. No. 928, vii, 95--113.
\newblock S{\'e}minaire Bourbaki. Vol. 2003/2004.

\bibitem{MR2507637}
{\sc Milman, E.}
\newblock On the role of convexity in isoperimetry, spectral gap and
  concentration.
\newblock {\em Invent. Math. 177}, 1 (2009), 1--43.

\bibitem{MR1008717}
{\sc Milman, V.~D., and Pajor, A.}
\newblock Isotropic position and inertia ellipsoids and zonoids of the unit
  ball of a normed {$n$}-dimensional space.
\newblock In {\em Geometric aspects of functional analysis (1987--88)},
  vol.~1376 of {\em Lecture Notes in Math.} Springer, Berlin, 1989,
  pp.~64--104.

\bibitem{MR856576}
{\sc Milman, V.~D., and Schechtman, G.}
\newblock {\em Asymptotic theory of finite-dimensional normed spaces},
  vol.~1200 of {\em Lecture Notes in Mathematics}.
\newblock Springer-Verlag, Berlin, 1986.
\newblock With an appendix by M. Gromov.

\bibitem{MR2276533}
{\sc Paouris, G.}
\newblock Concentration of mass on convex bodies.
\newblock {\em Geom. Funct. Anal. 16}, 5 (2006), 1021--1049.

\bibitem{MR2833584}
{\sc Paouris, G.}
\newblock Small ball probability estimates for log-concave measures.
\newblock {\em Trans. Amer. Math. Soc. 364}, 1 (2012), 287--308.

\bibitem{MR0315079}
{\sc Pr{\'e}kopa, A.}
\newblock Logarithmic concave measures with application to stochastic
  programming.
\newblock {\em Acta Sci. Math. (Szeged) 32\/} (1971), 301--316.

\bibitem{MR0404557}
{\sc Pr{\'e}kopa, A.}
\newblock On logarithmic concave measures and functions.
\newblock {\em Acta Sci. Math. (Szeged) 34\/} (1973), 335--343.

\bibitem{MR1694526}
{\sc Rudelson, M.}
\newblock Random vectors in the isotropic position.
\newblock {\em J. Funct. Anal. 164}, 1 (1999), 60--72.

\bibitem{MR2446328}
{\sc Sodin, S.}
\newblock An isoperimetric inequality on the {$l_p$} balls.
\newblock {\em Ann. Inst. Henri Poincar\'e Probab. Stat. 44}, 2 (2008),
  362--373.

\bibitem{SV}
{\sc Srivastava, N., and Vershynin, R.}
\newblock Covariance estimation for distributions with $2 + \varepsilon$
  moments.
\newblock {\em Ann. Prob.\/} (to appear).

\bibitem{MR517198}
{\sc Sudakov, V.~N.}
\newblock Typical distributions of linear functionals in finite-dimensional
  spaces of high dimension.
\newblock {\em Dokl. Akad. Nauk SSSR 243}, 6 (1978), 1402--1405.

\bibitem{MR0365680}
{\sc Sudakov, V.~N., and Cirel{\cprime}son, B.~S.}
\newblock Extremal properties of half-spaces for spherically invariant
  measures.
\newblock {\em Zap. Nau\v cn. Sem. Leningrad. Otdel. Mat. Inst. Steklov. (LOMI)
  41\/} (1974), 14--24, 165.
\newblock Problems in the theory of probability distributions, II.

\bibitem{MR2178341}
{\sc Vempala, S.}
\newblock Geometric random walks: a survey.
\newblock In {\em Combinatorial and computational geometry}, vol.~52 of {\em
  Math. Sci. Res. Inst. Publ.} Cambridge Univ. Press, Cambridge, 2005,
  pp.~577--616.

\bibitem{MR2853824}
{\sc Vempala, S.~S.}
\newblock Recent progress and open problems in algorithmic convex geometry.
\newblock In {\em 30th {I}nternational {C}onference on {F}oundations of
  {S}oftware {T}echnology and {T}heoretical {C}omputer {S}cience}, vol.~8 of
  {\em LIPIcs. Leibniz Int. Proc. Inform.} Schloss Dagstuhl. Leibniz-Zent.
  Inform., Wadern, 2010, pp.~42--64.

\bibitem{MR1440135}
{\sc von Weizs{\"a}cker, H.}
\newblock Sudakov's typical marginals, random linear functionals and a
  conditional central limit theorem.
\newblock {\em Probab. Theory Related Fields 107}, 3 (1997), 313--324.

\bibitem{Youssef}
{\sc Youssef, P.}
\newblock Estimating the covariance of random matrices.
\newblock {\em http://arxiv.org/abs/1301.6607\/} (preprint).

\end{thebibliography}
}

\address
\end{document}